\documentclass{elsarticle}
\usepackage{lineno,hyperref}
\modulolinenumbers[5]
\usepackage{booktabs,makecell}
\usepackage{multicol}
\usepackage{multirow}
\usepackage{amsmath,stmaryrd}
\usepackage{graphicx}
\usepackage{subfigure}
\usepackage{bm}
\usepackage{enumerate}
\usepackage{pgfplots}
\usepackage{tikz}
\usepackage{algorithm}
\usepackage[noend]{algpseudocode}
\biboptions{numbers,sort&compress}
\usepackage{color,xcolor,colortbl}
\usepackage{footnote}
\usepackage{url}
\usepackage{array}
\usepackage{etex}

\begin{document}
	
\begin{frontmatter}

		\title{Some Conclusions on Markov Matrices and Transformations}
		
        \author[rvt]{Chengshen Xu\corref{cor1}}
        \cortext[cor1]{Corresponding author: xcssgzs@126.com}
        \address[rvt]{Autohome Inc., 10th Floor Tower B, No. 3 Dan Ling Street Haidian District, Beijing, China}

\begin{abstract}

Markov matrices have an important role in the filed of stochastic processes. In this paper, we will show and prove a series of conclusions on Markov matrices and transformations rather than pay attention to stochastic processes although these conclusions are useful for studying stochastic processes. These conclusions we come to, which will make us have a deeper understanding of Markov matrices and transformations, refer to eigenvalues, eigenvectors and the structure of invariant subspaces. At the same time, we account for the corresponding significances of the conclusions. For any Markov matrix and the corresponding transformation, we decompose the space as a direct sum of an eigenvector and an invariant subspace. Enlightened by this, we achieve two theorems about Markov matrices and transformations inspired by which we conclude that Markov transformations may be a defective matrix--in other words, may be a nondiagonalizable one. Specifically, we construct a nondiagonalizable Markov matrix to exhibit our train of thought.
			
\end{abstract}
		
		
\end{frontmatter}
	
\section{Introduction}

Markov matrices have a widespread use in the filed of stochastic processes which are import in statistics \cite{heller1965stochastic}. As is well known, a matrix
\begin{eqnarray}
\label{Markov}
A=
\left(
  \begin{array}{cccc}
    a_{11} & a_{12} & \cdots & a_{1n} \\
    a_{21} & a_{22} & \cdots & a_{2n} \\
    \vdots & \vdots &   & \vdots \\
    a_{n1} & a_{n2} & \cdots & a_{nn} \\
  \end{array}
\right)
\end{eqnarray}
is called a Markov one if $\forall 1\leq i, j\leq n$ and the conditions \cite{kemeny1976markov, seneta2006non}
\begin{eqnarray}
\label{condition1}
&&
a_{ij}\geq 0,\\
\label{condition2}
&&
\sum\limits_{i=1}^{n}a_{ij}=1
\end{eqnarray}
are satisfied. In the science of probability and statistics, the matrix element $a_{ij}$ for a Markov matrix denotes the transition probability from state $j$ to state $i$, which is nonnegative as a matter of course (Eq. \ref{condition1}). The sum of the transition probabilities from state $j$ to the other states is $1$, as is expectational (Eq. \ref{condition2}).

In Section 2, we will study the Markov matrices and transformations themselves which belong to linear algebra rather than the stochastic processes which belong to probability and statistics. We will propose, demonstrate and discuss a series of conclusions that involve eigenvalues, eigenvectors, invariant subspaces and diagonalization about Markov matrices. At the same time, we will account for the corresponding significance of the conclusions. Firstly, we will demonstrate that any Markov matrix (or transformation) has at least one eigenvector whose any component is nonnegative. As is known, this eigenvector denotes the steady state of a stochastic process. Secondly, we will demonstrate that any Markov matrix (or transformation) has an invariant subspace, which is a hyperplane. Thirdly, we will demonstrate that for any Markov transformation, the space can be denoted as a direct sum of the eigenvector and the invariant subspace above. Finally, we will demonstrate two necessary and sufficient conditions of that a matrix is a Markov matrix.

In Section 3, we will discuss the significance of the eigenvector whose all components are nonnegative and the problem of diagonalization of Markov matrices. In this step, we will construct a nondiagonalizable Markov matrix to demonstrate the conclusion that a Markov matrix may be defective and to exhibit our train of thought. All the above will make us understand eigenvalues and eigenvectors of Markov matrices and the structures of Markov transformations more deeply.

\section{Theorems and Demonstrations}

Supposing matrix $A$ (Eq. \ref{Markov}) is a Markov matrix and $\lambda$ is an eigenvalue of it, because of Gershgorin circle theorem, for $A^{T}$ (the transposition matrix of $A$) we have eigenvalue \cite{weisstein2003gershgorin, bordenave2012circular}
\begin{eqnarray}
\lambda\in\bigcup\limits_{i=1}^{n}\{x||x-a_{ii}|\leq\sum\limits_{j=1}^{n}a_{ij}\}\Rightarrow -1\leq\lambda\leq 1.
\end{eqnarray}
which is also an eigenvalue of $A$. We write down the space as $V$ and the dimension of $V$ as $n$ and define
\begin{eqnarray}
&&
S1_{n}=\{X=(x_{1}, x_{2}, \cdots, x_{n})^{T}|\sum\limits_{i=1}^{n}x_{i}=1\},\\
&&
S2_{n}=\{X=(x_{1}, x_{2}, \cdots, x_{n})^{T}|\sum\limits_{i=1}^{n}x_{i}=1, x_{i}\geq 0, \forall 1\leq i\leq n\},\\
&&
S3_{n}=\{X=(x_{1}, x_{2}, \cdots, x_{n})^{T}|\sum\limits_{i=1}^{n}x_{i}=0\},\\
&&
S4_{n}=\{X=(x_{1}, x_{2}, \cdots, x_{n})^{T}|\sum\limits_{i=1}^{n}x_{i}=0, x_{i}\geq -\frac{1}{n}, \forall 1\leq i\leq n\}.
\end{eqnarray}
Evidently, $S4$ is the vertical projection of $S2_{n}$ at the hyperplane $S3$.

Now we propose a series of theorems as follows:

Theorem 1. $\forall$ Markov transformation $A$ and vector $X$, the transformation does not change the sum of the components of the vector. Specially, if $X\in S1_{n}$ is satisfied, we have $AX\in S1_{n}$.

Theorem 2. If vector $X\in S2_{n}$ is satisfied and $A$ is a Markov matrix, we have $AX\in S2_{n}$.

Theorem 3. Markov matrix $A$ has at least one eigenvector whose components are all nonnegative and the sum of whose components is positive. We write this eigenvector as $\bm{\alpha}$ and we set
\begin{eqnarray}
\bm{\xi}=\frac{\bm{\alpha}}{\sum\limits_{i=1}^{n}\alpha_{i}}\in S2_{n}
\end{eqnarray}

Theorem 4. Supposing $A$ is a Markov matrix and $\bm{\varepsilon}$ is an eigenvector of $A$, if $\sum\limits_{i=1}^{n}\varepsilon_{i}\neq 0$ is satisfied, we have the corresponding eigenvalue
\begin{eqnarray}
\lambda=1
\end{eqnarray}

Theorem 5. $S3$ is an invariant subspace of Markov matrices (or Markov transformations).

Theorem 6. $V$ can be decompose as a direct sum of $\bm{\xi}$ and S3, i.e.
\begin{eqnarray}
V=\bm{\xi}\bigoplus S3
\end{eqnarray}

Theorem 7. The necessary and sufficient condition for that A is a Markov matrix is that $\forall X\in S2_{n}$, $AX\in S2_{n}$ is satisfied.

Theorem 8. we define
\begin{eqnarray}
S5_{n}=\{X|X=Y-\bm{\xi}, Y\in S2_{n}\}
\end{eqnarray}
evidently, $S5_{n}\subseteq S3_{n}$
The necessary and sufficient condition for that A is a Markov matrix is that there is an eigenvector $\bm{\xi}$ of $A$ as described in Theorem 3, and $\forall X\in S5_{n}$, $AX\in S5_{n}$ is satisfied.

Now we demonstrate the theorems in last section.

1). Supposing $A$ is an Markov transformation, because $AX$ is equivalent to the linear combination of column vectors of $A$ by the components of $X$ and $\sum\limits^{n}_{i=1}a_{ij}=1$ is satisfied, the sum of the components of the transformed vector is not changed. Specifically when we set $Y=AX$, we have
\begin{eqnarray}
\sum\limits_{i=1}^{n}y_{i}=\sum\limits_{i=1}^{n}\sum\limits_{j=1}^{n}a_{ij}x_{j}=\sum\limits_{j=1}^{n}\sum\limits_{i=1}^{n}a_{ij}x_{j}=\sum\limits_{j=1}^{n}x_{j}.
\end{eqnarray}
Specially, if $\sum\limits_{j=1}^{n}x_{j}=1$ is satisfied, we have $\sum\limits_{i=1}^{n}y_{i}=1$, i.e. if $X\in S1_{n}$ is satisfied, we have $AX\in S1_{n}$

2). We set $Y=AX$. Because $A$ is a Markov matrix and $X\in S2_{n}$ is satisfied, we know $X\in S1_{n}$. From Theorem 1 we have $Y\in X1$. Because $\forall i, j$, $a_{ij}\geq 0$ and $x_{j}\geq 0$ are known, we have
\begin{eqnarray}
y_{i}=\sum\limits_{i=1}^{n}a_{ij}x_{j}\geq 0.
\end{eqnarray}
To sum up,
\begin{eqnarray}
Y=AX\in S2_{n}
\end{eqnarray}
is satisfied.

3). Because $S2_{n}$ is an $(n-1)$-dimensional compact set whose topological structure is the same as
\begin{eqnarray}
\{Z=(z_{1}, z_{2}, \cdots, z_{n-1})^{T}|0\leq z_{i}\leq 1, \forall 1\leq i\leq n-1\},
\end{eqnarray}
via an appropriate transformation of coordinates
\begin{eqnarray}
\label{coordinatetransformation}
\left\{
  \begin{array}{ll}
    x_{1}=f_{1}(z_{1}, z_{2}, \cdots, z_{n-1}) \\
    x_{2}=f_{2}(z_{1}, z_{2}, \cdots, z_{n-1}) \\
    ~~~~~~~~~~~~~~\vdots, \\
    x_{n}=f_{n}(z_{1}, z_{2}, \cdots, z_{n-1}) \\
  \end{array}
\right.
\end{eqnarray}
in $S2_{n}$, we can make
\begin{eqnarray}
S2_{n}=\{Z=(z_{1}, z_{2}, \cdots, z_{n-1})^{T}|0\leq z_{i}\leq 1, \forall 1\leq i\leq n-1\}
\end{eqnarray}
and $Y=F(Q)=AX=AF(Z)$, i.e. $Q=F^{-1}AF(Z)$. We set $G=F^{-1}AF$, specifically, so we have
\begin{eqnarray}
\label{daizheng_orig}
\left\{
  \begin{array}{ll}
    q_{1}=g_{1}(z_{1}, z_{2}, \cdots, z_{n-1}) \\
    q_{2}=g_{2}(z_{1}, z_{2}, \cdots, z_{n-1}) \\
    ~~~~~~~~~~~~~~\vdots, \\
    q_{n-1}=g_{n-1}(z_{1}, z_{2}, \cdots, z_{n-1}) \\
  \end{array}
\right..
\end{eqnarray}
Thus demonstrating Theorem 3 that Markov matrix $A$ has an eigenvector whose all componets are nonnegative is equivalent to demonstrating that there is at least one fixed point of $G(Z)$ in $S2_{n}$.
We use mathematical induction to demonstrate this theorem. When $n=2$, Eq. (\ref{daizheng_orig}) is
\begin{eqnarray}
q_{1}=g_{1}(z_{1}),
\end{eqnarray}
so there is evidently at least a fixed point of it in $S2_{2}=\{Z=(z_{1})|0\leq z_{1}\leq 1\}$. Now we assume there is also a fixed point of the Markov transformation in $S2_{n}$ for $n=k$. For $n=k+1$, Eq. (\ref{daizheng_orig}) becomes
\begin{eqnarray}
\label{daizheng}
\left\{
  \begin{array}{ll}
    q_{1}=g_{1}(z_{1}, z_{2}, \cdots, z_{k}) \\
    q_{2}=g_{2}(z_{1}, z_{2}, \cdots, z_{k}) \\
    ~~~~~~~~~~~~~~\vdots, \\
    q_{k}=g_{k}(z_{1}, z_{2}, \cdots, z_{k}) \\
  \end{array}
\right.
\end{eqnarray}
Evidently, $\forall$ fixed $z_{1}, z_{2}, \cdots, z_{k-1}$, $q_{k}=g_{k}(z_{k})$ have at least one fixed point. At this point, there is no harm to suppose
\begin{eqnarray}
\label{fixedplane}
z_{k}=h(z_{1}, z_{2}, \cdots, z_{k-1}).
\end{eqnarray}
It is a $(k-1)$-dimensional hyperplane whose topological structure is the same as $S2_{k}$. By substituting Eq. (\ref{fixedplane}) into the first $k-1$ equations of Eqs. (\ref{daizheng}) we can obtain $Q=I(Z)$, specifically
\begin{eqnarray}
\label{daizheng2}
\left\{
  \begin{array}{ll}
    q_{1}=i_{1}(z_{1}, z_{2}, \cdots, z_{k-1}) \\
    q_{2}=i_{2}(z_{1}, z_{2}, \cdots, z_{k-1}) \\
    ~~~~~~~~~~~~~~\vdots, \\
    q_{k-1}=i_{k-1}(z_{1}, z_{2}, \cdots, z_{k-1}) \\
  \end{array}
\right..
\end{eqnarray}
By induction hypothesis, we know there is at least one fixed point of Eq. (\ref{daizheng2})
\begin{eqnarray}
\left\{
  \begin{array}{ll}
    q_{1}=z_{1} \\
    q_{2}=z_{2} \\
    ~~~\vdots \\
    q_{k-1}=z_{k-1} \\
  \end{array}
\right..
\end{eqnarray}
So combining with $q_{k}=z_{k}$, there is at least one fixed point of Eq. (\ref{daizheng}), which implies the conclusion of Theorem 3 is also true for $n=k+1$.

4). if $A$ is a Markov matrix and $\bm{\varepsilon}$ is an eigenvector of $A$, supposing the corresponding eigenvalue is $\lambda$, from Theorem 1, we have
\begin{eqnarray}
\sum\limits_{i=1}^{n}\lambda\varepsilon_{i}=\sum\limits_{i=1}^{n}\varepsilon_{i}\Leftrightarrow \lambda\sum\limits_{i=1}^{n}\varepsilon_{i}=\sum\limits_{i=1}^{n}\varepsilon_{i}.
\end{eqnarray}
Because of $\sum\limits_{i=1}^{n}\varepsilon_{i}\neq 0$, we have $\lambda=1$.

5). Supposing that A is a Markov matrix, from Theorem 1, if $X\in S3_{n}$ is satisfied, we have $AX\in S3_{n}$. $\forall X, Y\in S3$, assuming $Z=k_{1}X+k_{2}Y$, we have
\begin{eqnarray}
\sum\limits_{i=1}^{n}z_{i}=\sum\limits_{i=1}^{n}(k_{1}x_{i}+k_{2}y_{i})=\sum\limits_{i=1}^{n}k_{1}x_{i}+\sum\limits_{i=1}^{n}k_{2}y_{i}
=k_{1}\sum\limits_{i=1}^{n}x_{i}+k_{2}\sum\limits_{i=1}^{n}y_{i}=0,
\end{eqnarray}
so $S3$ is closed in terms of addition and scalar-multiplication. To sum up, $S3_{n}$ is an invariant subspace of $A$.

6). Because the dimension of the invariant subspace $S3_{n}$ is equal to $n-1$ and from Theorem 3 there is an eigenvector $\bm{\xi}$ whose all components are nonnegative, we know $\bm{\xi}\notin S3_{n}$. So the space $V$ can be decomposed as the direct sum of $\bm{\xi}$ and $S3_{n}$, i.e.
\begin{eqnarray}
\label{directsum}
V=\bm{\xi}\bigoplus S3_{n}
\end{eqnarray}

7). The necessity is evident, which is Theorem 2, so we only need to demonstrate the sufficiency. To demonstrate $A$ is a Markov matrix, we only need to demonstrate that $\forall 1\leq i, j\leq n$, we have $a_{ij}>0$ on the one hand and $\forall 1\leq j\leq n$ we have $\sum\limits_{i=1}^{n}a_{ij}=1$ on the other hand. we use reduction to absurdity. If $\exists 1\leq i, j\leq n$, $a_{ij}<0$ is satisfied, for $X=(0, 0, \cdots, 0, 1, 0, 0, \cdots, 0)^{T}\in S2_{n}$ whose $j$-th component is $1$ and all the other components are $0$, we have
\begin{eqnarray}
AX=(a_{1j}, a_{2j}, \cdots, a_{ij}, \cdots, a_{nj})^{T}\notin S2_{n}
\end{eqnarray}
because of $a_{ij}<0$. This is in contradiction with $AX\in S2_{n}$, so $\forall 1\leq i, j\leq n$, $a_{ij}\geq 0$ is satisfied. If $\exists 1\leq j\leq n$, $\sum\limits_{i=1}^{n}a_{ij}\neq 1$ is satisfied, for $X=(0, 0, \cdots, 0, 1, 0, 0, \cdots, 0)^{T}\in S2_{n}$ whose $j$-th component is $1$ and all the other components are $0$, we have
\begin{eqnarray}
AX=(a_{1j}, a_{2j}, \cdots, a_{ij}, \cdots, a_{nj})^{T}\Rightarrow AX\notin S2_{n}
\end{eqnarray}
because of $\sum\limits_{i=1}^{n}a_{ij}\neq 1$. This is in contradiction with $AX\in S2_{n}$, so $\forall 1\leq j\leq n$, $\sum\limits_{i=1}^{n}a_{ij}=1$ is satisfied. To sum up, $A$ is a Markov matrix, so the sufficiency is satisfied.

8). Because of Eq. (\ref{directsum}), $\forall \bm{\alpha}\in S2_{n}$, it is obvious that we can decompose it as
\begin{eqnarray}
\bm{\alpha}=\bm{\beta}+\bm{\xi},
\end{eqnarray}
where $\bm{\beta}\in S5_{n}$. We have
\begin{eqnarray}
A\bm{\alpha}=A\bm{\beta}+A\bm{\xi}=A\bm{\beta}+\bm{\xi}
\end{eqnarray}
It is obvious that the necessary and sufficient condition for $A\bm{\alpha}\in S2_{n}$ is $A\bm{\beta}\in S5_{n}$, so from Theorem 7 we know the conclusion of Theorem 8 is true.

\section{Discussions of Two Attracting Problems}

At first, we discuss the eigenvector $\bm{\xi}$ of Markov matrix $A$. $\forall \bm{\alpha}\in V$, it is obvious that we can decompose it as
\begin{eqnarray}
\bm{\alpha}=\bm{\beta}+k\bm{\xi},
\end{eqnarray}
where $\bm{\beta}\in S3_{n}$. If $\forall 1\leq i\leq n$, $\alpha_{i}\geq 0$ and $\sum\limits_{i=1}^{n}\alpha_{i}=1$ are satisfied, $\bm{\alpha}$ can denote the probability distribution in a stochastic process and we have $k>0$ evidently. Supposing we have select a appropriate basis vector group $\bm{\eta}_{1}, \bm{\eta}_{2}, \cdots, \bm{\eta}_{n-1}$ in the invariant subspace $S3_{n}$ that make the transformation matrix $A$ is Jordan standard form one \cite{browne1940reduction}
\begin{eqnarray}
A=
\left(
  \begin{array}{ccccc}
    1 & &  &  &  \\
    & J_{n_{1}}(\lambda_{1}) &  &  &  \\
     & & J_{n_{2}}(\lambda_{2}) &  &  \\
     & &  & \ddots &  \\
     & &  &  & J_{n_{s}}(\lambda_{s}) \\
  \end{array}
\right).
\end{eqnarray}
under the basis vector group $\bm{\xi}, \bm{\eta}_{1}, \bm{\eta}_{2}, \cdots, \bm{\eta}_{n-1}$.
Here
\begin{eqnarray}
&&
J_{n_{i}}(\lambda_{i})=
\left(
  \begin{array}{ccccc}
    \lambda_{i} & 1 & 0 & \cdots & 0 \\
    0 & \lambda_{i} & 1 & \cdots & 0 \\
    0 & 0 & \lambda_{i} & \cdots & 0 \\
    \vdots & \vdots & \vdots &  & \vdots \\
    0 & 0 & 0 & \cdots & \lambda_{i} \\
  \end{array}
\right).
\\
&&
\sum\limits_{i=1}^{s}n_{i}+1=n
\end{eqnarray}
is satisfied, where $n_{i}$ is the dimension of $J_{n_{i}}(\lambda_{i})$. Assuming under the basis vector group $\bm{\xi}, \bm{\eta}_{1}, \bm{\eta}_{2}, \cdots, \bm{\eta}_{n-1}$
\begin{eqnarray}
\bm{\beta}=(\beta_{1}, \beta_{2}, \cdots, \beta_{n})^{T}
\end{eqnarray},
we have $\beta_{1}=0$ because of $\bm{\beta}\in S3_{n}=W$. One can transform $\bm{\alpha}$ repeatedly via $A$,
\begin{eqnarray}
&&
\lim\limits_{m\rightarrow\infty}A^{m}\bm{\alpha}=\lim\limits_{m\rightarrow\infty}A^{m}\bm{\beta}+\lim\limits_{m\rightarrow\infty}A^{m}k\bm{\xi}
\nonumber\\
&&
=
\lim\limits_{m\rightarrow\infty}
\left(
  \begin{array}{ccccc}
    1 & &  &  &  \\
    & J_{n_{1}}^{m}(\lambda_{1}) &  &  &  \\
     & & J_{n_{2}}^{m}(\lambda_{2}) &  &  \\
     & &  & \ddots &  \\
     & &  &  & J_{n_{s}}^{m}(\lambda_{s}) \\
  \end{array}
\right)
\left(
  \begin{array}{c}
    0 \\
    \beta_{2} \\
    \vdots \\
    \beta_{n} \\
  \end{array}
\right)
+k\bm{\xi}.
\end{eqnarray}
In a stochastic process, a vector denotes the state distribution. For $-1<\lambda_{i}< 1$,
\begin{eqnarray}
\lim\limits_{m\rightarrow\infty}J_{n_{i}}^{m}(\lambda_{i})=O_{n_{i}}
\end{eqnarray}
where $O_{n_{i}}$ is $n_{i}\times n_{i}$ zero matrix. For $\lambda_{i}=-1$, $\lim\limits_{m\rightarrow\infty}J_{n_{i}}^{m}(-1)$ is diverged. In this situation if $\exists \sum\limits_{l=0}^{i-1}n_{l}<j\leq \sum\limits_{l=0}^{i}n_{l}$ (assuming $n_{0}=1$ for convenience), $\beta_{j}\neq 0$ is satisfied, the steady state cannot be reached. For $\lambda=1$ and $n_{i}> 1$, $J_{n_{i}}(1)$ is also diverged. In this situation, if $\exists \sum\limits_{l=0}^{i-1}n_{l}+1<j\leq \sum\limits_{l=0}^{i}n_{l}$, $\beta_{j}\neq 0$ is satisfied, the steady state cannot also be reached. In the situation the steady state can be reached, if $\forall i, \lambda_{i}\neq 1$, the final steady state is $k\bm{\xi}$ evidently, but if some $\lambda_{i}=1$, the final steady state is determined by the initial vector $\bm{\alpha}$. Because the dimension of the eigensubspace of $\lambda=1$ is greater than $1$, the steady state is not only in this situation.

One may also doubt the problem that whether any Markov matrix $A$ can be diagonalized or not, because studying it can make us have a deeper understanding of the structures of Markov matrices. In Theorem 6, we have decompose the space as the direct sum of $\bm{\xi}$ and $S3_{n}$, so we only need to research the effects of Markov matrices on the invariant subspace $S3_{n}$. $\forall$ vector $\bm{\alpha}\in S2_{n}$, it is obvious that we can decompose it as
\begin{eqnarray}
\bm{\alpha}=\bm{\beta}+\bm{\xi},
\end{eqnarray}
where $\bm{\beta}\in S5_{n}\subseteq S3_{n}$. In order to construct a nondiagonalizable Markov matrix, the key point is to construct a nondiagonalizable transformation which garantees $A\bm{\beta}\in S5_{n}$ for $\bm{\beta}\in S5_{n}$ enlightened by Theorem 8. To show the train of thought, we construct a nondiagonalizable Markov matrix in 3-dimensional space. For simplicity,
we set
\begin{eqnarray}
\bm{\xi}=(\frac{1}{3}, \frac{1}{3}, \frac{1}{3})^{T},
\end{eqnarray}
so we have
\begin{eqnarray}
S5_{3}=S4_{3}=\{X=(x_{1}, x_{2}, x_{3})^{T}|\sum\limits_{i=1}^{3}x_{i}=0, x_{i}\geq -\frac{1}{3}, \forall 1\leq i\leq 3\},
\end{eqnarray}
which is the vertical projection of $S2_{3}$ on $S3_{3}$. In order to construct a Markov matrix (or transformation) we set the basis vector group
\begin{eqnarray}
\left\{
  \begin{array}{ll}
\bm{\varepsilon}_{1}=\bm{\xi}=(\frac{1}{3}, \frac{1}{3}, \frac{1}{3})^{T},\\
\bm{\varepsilon}_{2}=(\frac{1}{2}, -\frac{1}{4}, -\frac{1}{4})^{T}\in S4_{3}\subseteq S3_{3},\\
\bm{\varepsilon}_{3}=(-\frac{1}{2}, 1 , -\frac{1}{2})^{T}\in S4_{3}\subseteq S3_{3},
  \end{array}
\right.
\end{eqnarray}
so the matrix constituted by the basis vector group is
\begin{eqnarray}
\left(
\begin{array}{ccc}
\frac{1}{3}  & \frac{1}{2}  & -\frac{1}{2} \\
\frac{1}{3}  & -\frac{1}{4} & 1 \\
\frac{1}{3} & -\frac{1}{4} & -\frac{1}{2}\\
\end{array}
\right)
\end{eqnarray}
and its inverse matrix is
\begin{eqnarray}
\left(
\begin{array}{ccc}
1  & 1   & 1  \\
\frac{4}{3}  & 0 & -\frac{4}{3} \\
0 & \frac{2}{3}  & -\frac{2}{3}\\
\end{array}
\right).
\end{eqnarray}
As is known, any matrix, which cannot be diagonalized, can be converted to the Jordan canonical form by a similarity transformation \cite{browne1940reduction}. Thus we assume that the Jordan canonical form of $A$ is
\begin{eqnarray}
\left(
\begin{array}{ccc}
1  & 0  & 0 \\
0  & 0 & 1 \\
0 & 0  & 0\\
\end{array}
\right).
\end{eqnarray}
Thus, when the basis vector group is
\begin{eqnarray}
\left\{
  \begin{array}{ll}
\bm{\eta}_{1}=(1, 0, 0)^{T},\\
\bm{\eta}_{2}=(0, 1, 0)^{T},\\
\bm{\eta}_{3}=(0, 0 , 1)^{T},
  \end{array}
\right.
\end{eqnarray}
we can obtain the the matrix of the transformation
\begin{eqnarray}
&&
A=
\left(
\begin{array}{ccc}
\frac{1}{3}  & \frac{1}{2}  & -\frac{1}{2} \\
\frac{1}{3}  & -\frac{1}{4} & 1 \\
\frac{1}{3} & -\frac{1}{4} & -\frac{1}{2}\\
\end{array}
\right)
\left(
\begin{array}{ccc}
1  & 0  & 0 \\
0  & 0 & 1 \\
0 & 0  & 0\\
\end{array}
\right)
\left(
\begin{array}{ccc}
1  & 1   & 1  \\
\frac{4}{3}  & 0 & -\frac{4}{3} \\
0 & \frac{2}{3}  & -\frac{2}{3}\\
\end{array}
\right)
\nonumber\\
&&
=
\left(
\begin{array}{ccc}
\frac{1}{3}  & \frac{2}{3}   & 0  \\
\frac{1}{3}  & \frac{1}{6} & \frac{1}{2} \\
\frac{1}{3} & \frac{1}{6}  & \frac{1}{2}\\
\end{array}
\right),
\end{eqnarray}
which is a nondiagonalizable Markov matrix. Of course, when the eigenvalue in $S3$ is not equal to $0$, we can also construct a nondiagonalizable one as long as the second basis vector is small enough.

\section{The Summary and Prospect}

In this paper, we propose and demonstrate a series of conclusions on Markov matrices and transformations, including the state vector $\bm{\xi}$, the invariant subspace and the direct sum decomposition of the space, which makes us have a deeper understanding of Markov matrices and transformations. The proposition of the two necessary and sufficient conditions for that A matrix is a Markov matrix provide a train of thought to construct a nondiagonalizable Markov matrix for us. Although the conclusions and discussions in this paper are more about linear algebra rather than probability and statistics, we think they will be useful in the domain of probability and statistics and even in the domain of machine learning.

	
	
\section*{References}
    \bibliographystyle{elsart-num1}
	\bibliography{research}

\begin{thebibliography}{1}
\expandafter\ifx\csname url\endcsname\relax
  \def\url#1{\texttt{#1}}\fi
\expandafter\ifx\csname urlprefix\endcsname\relax\def\urlprefix{URL }\fi

\bibitem{heller1965stochastic}
A.~Heller, On stochastic processes derived from markov chains, The Annals of
  Mathematical Statistics 36~(4) (1965) 1286--1291.

\bibitem{kemeny1976markov}
J.~G. Kemeny and J.~L. Snell, Markov chains, Springer-Verlag, New York, 1976.

\bibitem{seneta2006non}
E.~Seneta, Non-negative matrices and Markov chains, Springer Science \&
  Business Media, 2006.

\bibitem{weisstein2003gershgorin}
E.~W. Weisstein, Gershgorin circle theorem.

\bibitem{bordenave2012circular}
C.~Bordenave, P.~Caputo and D.~Chafa{\"\i}, Circular law theorem for random
  markov matrices, Probability Theory and Related Fields 152~(3-4) (2012)
  751--779.

\bibitem{browne1940reduction}
E.~Browne, On the reduction of a matrix to a canonical form, The American
  Mathematical Monthly 47~(7) (1940) 437--450.

\end{thebibliography}
	
\end{document}